\documentclass[11pt,a4paper]{amsart}
\usepackage{graphics}

\usepackage{amsmath}

\usepackage{amssymb}

\theoremstyle{plain}   
\begingroup 

\newtheorem{theorem}{Theorem}[section]   
\newtheorem{corollary}[theorem]{Corollary}     
\newtheorem{lemma}[theorem]{Lemma}         
\newtheorem{proposition}[theorem]{Proposition}  
\endgroup

\theoremstyle{definition}

\theoremstyle{remark}
\newtheorem{remark}[theorem]{Remark}        
\newtheorem{example}[theorem]{Example}        

\numberwithin{equation}{section}

\newcommand{\sign}{\operatorname{sign}}

\newcommand{\R}{{\mathbb R}}
\newcommand{\N}{{\mathbb N}}

\newcommand{\e}{\varepsilon}

\newcommand{\vp}{\varphi}

\newcommand{\diam}{\operatorname{diam}}

\newcommand{\Hau}{{\mathcal{H}}}

\newcommand{\mycon}{q}
\newcommand{\complex}{\mathbb{C}}
\newcommand{\h}{h}
\newcommand{\dd}{\rho}
\newcommand{\gap}{\theta}

\usepackage[dvips]{color}

\begin{document}
\title[Metric derived numbers]
{Metric derived numbers and continuous metric differentiability via
homeomorphisms}

\author{Jakub Duda}
\author{Olga Maleva}

\thanks{This paper is dedicated to Nigel Kalton on the occasion of
his sixtieth birthday.}



\email{duda@karlin.mff.cuni.cz}
\thanks{The first author was supported in part by ISF}
\email{o.maleva@dpmms.cam.ac.uk}

\address{
Charles University\\
Department of Mathematical Analysis\\
Sokolovsk\'a 83\\
186 75 Praha 8\\
Czech Republic
}

\address{
Department of Pure Mathematics and Mathematical Statistics\\
Centre for Mathematical Sciences\\ 
University of Cambridge\\ 
Wilberforce Road\\ 
Cambridge\\ 
CB3 0WB \\
United Kingdom
}

%
\date{August 14, 2006}
%
%
\begin{abstract} 
We define the notions of unilateral metric derivatives and
``metric derived numbers'' in analogy
with Dini derivatives (also referred to as ``derived numbers'') 
and establish their basic properties.
We also prove that the set of points where a path with values
in a metric space with continuous metric derivative is not
``metrically differentiable'' (in a certain strong sense)
is $\sigma$-symmetrically porous and provide an example
of a path for which this set is uncountable.
In the second part of this paper, 
we study the continuous metric differentiability via a homeomorphic
change of variable.
\end{abstract}

\maketitle

\section{Introduction}
The main aim of this paper is to study analogues of the usual notion of
differentiability which work for mappings with values in metric
spaces.
Let $(X,\rho)$ be a metric space and $f:[a,b]\to X$ be any mapping. As
every metric space isometrically embeds in some Banach space (see
e.g.\ \cite[Lemma~1.1]{BL}), we can suppose that the distance in $X$
is in fact generated by a complete norm $\|\cdot\|$.
Define 
\[ md_\pm(f,x)=\lim_{t\to0+} \frac{\|f(x\pm t)-f(x)\|}{t}\]
to be the {\em unilateral right $($resp.\ left$)$ metric derivatives}
of the mapping $f$ at $x$. 
If $md_+(f,x)$ and $md_-(f,x)$
exist, and are equal, then we call $md(f,x):=md_+(f,x)$ the metric
derivative of $f$ at the point $x$.

We say that {\em $f$ is metrically differentiable
at~$x$} provided $md(f,x)$ exists and 
\begin{equation}\label{md}
\|f(y)-f(z)\|-md(f,x)|y-z|=o(|y-x|+|z-x|),\text{ when }(y,z)\to(x,x).
\end{equation}
Note that in this terminology, the 
existence of the ``metric derivative'' $md(f,x)$ of $f$ at $x$ does 
not necessarily imply that $f$ is
metrically differentiable at $x$! The basic example of such mapping
would be $f(t)=|t|:\R\to\R$ and $x=0$.

Metric derivatives were introduced by Kirchheim in~\cite{K}
(see also~\cite{Am,DP,KS}),
and were studied by several authors (see e.g.~\cite{AK,Dac,Dkir,DZ2}). 
In~\cite{AK}, the authors work with a slightly weaker version
of metric differentiability.

We start section~\ref{unilatsec} by noting
that the set of points where $md_\pm(f,x)$ exist, but $md_+(f,x)\neq
md_-(f,x)$, 
is countable; see Theorem~\ref{unilatthm}. This is analogous
to a similar theorem for unilateral derivatives of real-valued
functions.

There is a well established theory of derived numbers (or Dini
derivatives) of real-valued functions $f:\R\to\R$ (see e.g.\ \cite{Br}). In 
section~\ref{unilatsec},
we generalize theorems about relationships among the Dini derivatives
to the context of metric derived numbers $mD^\pm, mD_\pm$.

In Theorem~\ref{angularptsthm}, we prove that the set of 
``angular'' points of each $f:\R\to X$, i.e.\ points $x\in\R$
where either
$mD_-(f,x)>mD^+(f,x)$ or $mD_+(f,x)>mD^-(f,x)$, is countable.
Theorem~\ref{upperdifthm} (resp.\ Theorem~\ref{lowerdifthm} if $f$ is pointwise-Lipschitz)
shows that the sets of points $x\in\R$ where $mD^+(f,x)\neq mD^-(f,x)$
(resp.\ $mD_+(f,x)\neq mD_-(f,x)$) is $\sigma$-porous.
Theorem~\ref{butterflythm} (see also Corollary~\ref{butterflycor}) 
is a metric analogue of the so-called Denjoy-Young-Saks
theorem about Dini derivatives (see e.g.\ \cite[Theorem~4.4]{Br}).

In section~\ref{ptsmetrnondifsec},
we show that if $md(f,\cdot)$ is
a continuous function, then the set of points $x$, where $f$ is not metrically differentiable,
is $\sigma$-symmetrically porous (Theorem~\ref{ptmetrnondif}). 
In Theorem~\ref{difexample}, we show that this set is not necessarily 
countable.
This means that the properties of metric derivatives are
different from the properties of standard ones; in the latter case,
the set considered in section~\ref{ptsmetrnondifsec},
would necessarily be countable
(if say $md(f,\cdot)\equiv1$ for a real-valued $f$ then the 
standard unilateral derivatives of $f$ are equal to $\pm1$ at all points).
\par
In section~\ref{auxsec}, we discuss sufficient conditions
for a mapping to be metrically differentiable at a point.
This is closely related to the notion of bilateral metric regularity.

In a recent paper~\cite{DZ2}, L. Zaj\'\i\v cek 
together with the first author
characterized those mappings $f:[a,b]\to X$ that allow
a metrically differentiable 
(resp.\ boundedly metrically differentiable) parameterization.
In section~\ref{contsec}, we study the situation when $f$
allows a 
continuously metrically differentiable parameterization
(by this we mean that for a suitable homeomorphism $h$,
the composition $f\circ h$ is metrically differentiable and
its metric derivative is continuous),
or just a parameterization with continuous metric derivative;
see Theorems~\ref{mdifctsthm} and~\ref{ctsmdthm} for more details.

\section{Preliminaries}\label{preliminaries}

By $\lambda$ we 
denote the $1$-dimensional Lebesgue
measure on $\R$, and by $\Hau^1$ the 
$1$-dimensional Hausdorff measure.
In the following, $X$ is always a real Banach space.
\par
The following is a version of the Sard's theorem.
For a proof see e.g.~\cite[Lemma~2.2]{DZ2}.

\begin{lemma}\label{mdsardlem}
Let $f:[a,b]\to X$ be arbitrary.
Then $\Hau^1(f(\{x\in[a,b]:md(f,x)=0\}))=0$.
\end{lemma}

By $B(x,r)$, we denote the open ball in $X$ with center $x\in X$ and radius $r>0$.
Let $M\subset\R$, $x\in M$, and $R>0$. Then we define $\gamma(x,R,M)$
to be the supremum  of all $r>0$  for which there exists $z\in\R$
such that $B(z,r)\subset B(x,R)\setminus M$. 
Also, we define $S\gamma(x,R,M)$ to be the 
supremum of all $r>0$ for which there exists $z\in\R$ 
such that $B(z,r)\cup B(2x-z,r)\subset B(x,R)\setminus M$.
Further, we 
define the {\em upper porosity} of $M$ at $x$ as
\[ \overline{p}(M,x):=2\limsup_{R\to0+}\frac{\gamma(x,R,M)}{R},\]
and the {\em symmetric upper porosity} of $M$ at $x$ as
\[ S\overline{p}(M,x):=2\limsup_{R\to0+}\frac{S\gamma(x,R,M)}{R}.\]
We say that $M$ is {\em porous\footnote{In the terminology of~\cite{Zsur},
this corresponds to $M$ being ``an upper-porous set''.}} (resp.\
{\em symmetrically porous}) provided 
$\overline{p}(M,x)>0$ for all $x\in M$
(resp.\ $S\overline{p}(M,x)>0$ for all $x\in M$). 
We say that $N\subset\R$ is
{\em $\sigma$-porous} (resp. {\em $\sigma$-symmetrically porous}) 
provided it is a countable union of porous 
(resp. symmetrically porous) sets.
For more information about porous sets, see a recent survey~\cite{Zsur}.
\par
Let $f:[a,b]\to X$. Then we say that {\em $f$ has finite variation}
or that {\em $f$ is BV},
provided $\bigvee^b_a f<\infty$. 
(Recall that $\bigvee^b_a f =\sup_D \sum_{i=0}^{n(D)-1} \|f(x_i)- f(x_{i+1})\|$,
where the supremum is taken over all 
partitions $D = \{a = x_0 < x_1 < \dots < x_n = b\}$, 
of $[a,b]$ and $n(D) = \# D - 1$.)
We define $\bigvee^u_v f:=-\bigvee^v_u f$
for $a\leq u<v\leq b$. We will denote $v_f(x):=\bigvee_a^x f$ for $x\in[a,b]$.

We say that $f:\R\to X$ is {\em pointwise-Lipschitz} if 
$\limsup_{y\to x}\frac{\|f(x)-f(y)\|}{|x-y|}<\infty$
for every $x\in \R$.

A considerable part of the present article is devoted to metric analogues of
derived numbers (Dini derivatives).  Now, we give a definition of
metric derived numbers. 
Let  $f:\R\to X$. Define 
\[ mD^\pm(f,x)=\limsup_{t\to0+} \frac{\|f(x\pm t)-f(x)\|}{t},\]
and
\[ mD_\pm(f,x)=\liminf_{t\to0+} \frac{\|f(x\pm t)-f(x)\|}{t},\]
to be the {\em unilateral upper $($resp.\ lower$)$ metric derived numbers}
(we also allow the value $+\infty$).

Note that 
if all four metric derived numbers of
a mapping $f:\R\to X$ agree at a point $x$, then $md(f,x)$ exists,
but still $f$ is not necessarily metrically differentiable at $x$.

\section{Unilateral metric derivatives}\label{unilatsec}

It is well known that the set where the standard unilateral derivatives
of a real function of a real variable exist but are not equal is countable
(see e.g.\ \cite[Theorem~7.2]{Jef}).
The following theorem shows that it is also true for unilateral metric 
derivatives.

\begin{theorem}\label{unilatthm}
Let $f:\R\to X$. Then the set of points $x\in\R$ where $md_+(f,x)$, $md_-(f,x)$
exist but $md_+(f,x)\neq md_-(f,x)$, is countable.
\end{theorem}

\begin{proof} The proof is similar to the proof of~\cite[Theorem~7.2]{Jef}
and thus we omit it.
\end{proof}


It is well known that for a real function of a real variable the
set of angular points (i.e.\ points where $D_-f>D^+f$ or $D_+f>D^-f$;
$D^\pm f $, $D_\pm f$ are the standard derived numbers) is countable;
see e.g.\ \cite[Theorem~7.2]{Jef}. The following theorem shows
what happens for metric derived numbers.

\begin{theorem}\label{angularptsthm}
Let $f:\R\to X$. Then the set of points $x\in\R$ where either $mD_-(f,x)>mD^+(f,x)$
or $mD_+(f,x)>mD^-(f,x)$ is countable.
\end{theorem}

\begin{proof}
By symmetry, it is enough to prove that the set
${E}=\{x\in\R:mD_-(f,x)>mD^+(f,x)\}$
is countable.
Let $h<k$ be two positive rational numbers. For a positive integer $n$ let
$E_{hkn}$ be the set of points 
$x\in E$ for which
$\frac{\|f(\xi)-f(x)\|}{|\xi-x|}<h$ and
$\frac{\|f(\xi')-f(x)\|}{|\xi'-x|}>k$
whenever
$0<\xi-x<1/n$ and $0<x-\xi'<1/n$. Then $E_{hkn}\cap (x-1/n,x+1/n)=\{x\}$.
Suppose that is not true, and there is a point 
$x_1\in E_{hkn}\cap (x-1/n,x+1/n)$ such that $x_1\ne x$. 
Then assuming $x>x_1$, say, we get
$\frac{\|f(x_1)-f(x)\|}{|x_1-x|}<h$ and 
$\frac{\|f(x)-f(x_1)\|}{|x-x_1|}>k$,
a contradiction. Thus all points of $E_{hkn}$ are isolated, 
and $E_{hkn}$ is countable. Because 
$E\subset\bigcup_{h,k,n}E_{hkn}$, we obtain
the conclusion of the theorem.
\end{proof}

We have the following two theorems concerning the points where unilateral
lower and upper metric derivatives differ. In the proofs, we use similar ideas
as in~\cite[Theorem~1]{EH}.

\begin{theorem}\label{upperdifthm}
Let $X$ be a Banach space, and $f:\R\to X$ be arbitrary. Then the set 
\[ \{x\in\R:mD^+(f,x)\neq mD^-(f,x)\}\]
is $\sigma$-porous.
\end{theorem}

\begin{proof}
We will only prove that the set
\[ A=A_f=\{x\in\R:mD^-(f,x)<mD^+(f,x)\},\]
is $\sigma$-porous (and notice that $\{x\in\R:mD^-(f,x)>mD^+(f,x)\}$ is 
$\sigma$-porous as it is equal to $A_{f(-\cdot)}$).
To that end, it is enough to establish that
\[ A_{rs}=\{x\in A: mD^-(f,x)<r<s<mD^+(f,x)\},\]
is $\sigma$-porous for all $r<s$ pairs of positive rational numbers.
Define
\begin{equation*}
A_{rsn}=\bigg\{x\in A_{rs}:\frac{\|f(x)-f(y)\|}{|x-y|}<r\text{ for }y\in(x-1/n,x)\bigg\}.
\end{equation*}
We easily see that $A_{rs}=\bigcup_n A_{rsn}$. We will prove that
$A_{rsn}$ is $\frac{\delta-1}{\delta}$-porous, where $\delta=\min(2,(s+r)/2r)$.
Let $x\in A_{rsn}$. Then there exist $x_k\to x+$ such that
$ \frac{\|f(x)-f(x_k)\|}{|x-x_k|}>s$. 
Choose $k$ large enough such that $|x-x_k|<1/n$.
Define $w_k=x+\delta\,(x_k-x)$, and let $y\in[x_k,w_k]\cap A_{rsn}$.
Then
\begin{equation*}
\begin{split}
\|f(x)-f(y)\|&\geq 
\|f(x)-f(x_k)\|-\|f(x_k)-f(y)\| \\
& \geq s\,|x-x_k|-r\,|x_k-y| \geq s\,|x-x_k|-r\,|x_k-w_k| \\
& = s\,|x-x_k|-r\,(\delta-1)\,|x_k-x|\\
& = |x-x_k|\,(s-r\,(\delta-1)) = |w_k-x|\,\frac{(s-r\,(\delta-1))}{\delta}\\
& \geq r\,|x-y|,
\end{split}
\end{equation*}
by the choice of $\delta$ (we used that $w_k-x=\delta\,(x_k-x)$, and $w_k-x_k=(\delta-1)\,(x_k-x)$).
Thus $y\not\in A_{rsn}$, and $[x_k,w_k]\cap A_{rsn}=\emptyset$.
Finally, note that
$\frac{w_k-x_k}{w_k-x}=\frac{\delta-1}{\delta}>0$.
\end{proof}

\begin{theorem}\label{lowerdifthm}
Let $X$ be a Banach space, and $f:\R\to X$ be pointwise-Lipschitz. Then the set 
\[ \{x\in\R:mD_+(f,x)\neq mD_-(f,x)\}\]
is $\sigma$-porous.
\end{theorem}

\begin{proof} 
We will only prove that the set  
\[B=B_f=\{x\in\R:mD_-(f,x)<mD_+(f,x)\},\]
is $\sigma$-porous,
and notice that $\{x\in\R:mD_-(f,x)>mD_+(f,x)\}$ is 
$\sigma$-porous as it is equal to $B_{f(-\cdot)}$.
We will prove that $B_f$
is $\sigma$-porous for $f$ that is pointwise-Lipschitz.
To that end, it is enough to establish that
$B_{rs}=\{x\in B: mD_-(f,x)<r<s<mD_+(f,x)\}$,
is $\sigma$-porous
for all $r<s$ pairs of positive rational numbers.
For $n\in\N$, define
\begin{equation*}
\begin{split}
B_{rsn}=\bigg\{x\in B_{rs}: & \frac{\|f(x)-f(y)\|}{|x-y|}>s\text{ for }y\in(x,x+1/n),\\
& \text{ and }\frac{\|f(x)-f(z)\|}{|x-z|}<n
\text{ whenever }0<|z-x|<1/n\bigg\}.
\end{split}
\end{equation*}
Since $f$ is pointwise-Lipschitz,
we easily see that $B_{rs}=\bigcup_n B_{rsn}$. 
We will prove that
$B_{rsn}$ is $\frac{\delta-1}{\delta}$-porous, where $\delta=\min(\frac{s-r}{n}+1,2)$.
Let $x\in B_{rsn}$. Then there exist $x_k\to x-$ such that
$\frac{\|f(x)-f(x_k)\|}{|x-x_k|}<r$. 
Choose $k$ large enough such that $|x-x_k|<1/n$.
Define $w_k=x-\delta\,(x-x_k)$, and let $y\in[w_k,x_k]\cap B_{rsn}$.
Then
\begin{equation*}
\begin{split}
\|f(x)-f(y)\|&\leq 
\|f(x)-f(x_k)\|+\|f(x_k)-f(y)\| \\
& \leq r\,|x-x_k|+n\,|x_k-y|  \leq r\,|x-x_k|+n\,|x_k-w_k| \\
& = r\,|x-x_k|+n\,(\delta-1)\,|x_k-x| = |x-x_k|\,(r+n\,(\delta-1))\\
& \leq s\,|x-y|,
\end{split}
\end{equation*}
by the choice of $\delta$ (we used that $x-w_k=\delta\,(x-x_k)$, $x_k-w_k=(\delta-1)\,(x-x_k)$,
and $|x_k-w_k|<1/n$).
Thus $y\not\in {B}_{rsn}$, and $[x_k,w_k]\cap {B}_{rsn}=\emptyset$.
Finally, note that
$\frac{x_k-w_k}{x-w_k}=\frac{\delta-1}{\delta}>0$.
\end{proof}

The following theorem asserts that outside of a set of measure $0$,
the fact that $mD^+(f,x)<\infty$ already implies that $md(f,x)$ exists.

\begin{theorem}\label{butterflythm} 
Let $f:\R\to X$ be arbitrary. Then there exists a set $N$ with Lebesgue measure
zero such that 
\begin{quote} 
if $x\in\R\setminus N$ and $mD^+(f,x)<\infty$,
then $md(f,x)$ exists, 
and $md(f,x)=mD^+(f,x)$.
\end{quote}
\end{theorem}

\begin{proof}
Let $N_1$ be the set of points $x\in\R$ where $mD^-(f,x)\neq mD^+(f,x)$. 
Then, by Theorem~\ref{upperdifthm} $N_1$ is $\sigma$-porous. 
Therefore, by the Lebesgue density
theorem, its Lebesgue measure $\lambda(N_1)$ is zero.
Let 
\[A_n=\{x\in\R:\|f(x+h)-f(x)\|\leq nh\text{ for }0<h<1/n\}.\]
Let $A$ be the set of points $x$ such 
that $mD^+(f,x)<\infty$. Then $A=\bigcup_n A_n$. Let $A_{n,j}$ be
subsets of $A_n$, such that $A_n=\bigcup A_{n,j}$, and $\diam(A_{n,j})<1/n$.
Then $f|_{A_{n,j}}$ is $n$-Lipschitz, and thus,
by Kirszbraun theorem, see \cite{Ki}, it can be extended
to an $n$-Lipschitz function $f_{n,j}$ defined on the whole real line.
By~\cite[Theorem~2.7]{Dkir}, we obtain that $f_{n,j}$ is metrically
differentiable at all $x\in D_{n,j}$, where $\lambda(\R\setminus D_{n,j})=0$.
Let $E_{n,j}\subset D_{n,j}\cap A_{n,j}$ be 
the set of points
of density of $D_{n,j}\cap A_{n,j}$. 
By the Lebesgue density theorem we have that 
$\lambda(D_{n,j}\cap A_{n,j}\setminus E_{n,j})=0$.
We shall prove that $md(f,x)$ exists and is equal to $mD^+(f,x)$ at all points $x\in E_{n,j}$
for all $n,j\in\N$. 
This will conclude the proof, as 
the set $N=\bigcup_{n,j}(A_{n,j}\setminus E_{n,j})$
has Lebesgue measure $0$.
\par
To finish the proof, let $x\in E_{n,j}$. Fix $\e>0$. Find $\delta>0$ such that
$\lambda(E_{n,j}\cap(x,x+t))\geq \big(1-\frac{\e}{4n}\big)t$
for $0<t<\delta$, and 
$\left| \frac{\|f_{n,j}(x+t)-f_{n,j}(x)\|}{|t|}-md(f_{n,j},x)\right|\leq\e$,
whenever $0<|t|<\delta$. 
Thus for each $0<h<\delta$ there exists $y\in E_{n,j}\cap (x,x+h)$
such that $|y-(x+h)|\leq \frac{\e h}{2n}$. Now,
\begin{equation*}
\begin{split}
\|f(x+h)-f(x)\|&\leq\|f(y)-f(x)\|+\|f(x+h)-f(y)\|\\
&\leq(md(f_{n,j},x)+\e)(y-x)+\e h\\
&\leq (md(f_{n,j},x)+\e)h,
\end{split}
\end{equation*}
since $x$ and $y$ belong to $E_{n,j}\subset A_n$ and $y>x$.
On the other hand,
\begin{equation*}
\begin{split}
\|f(x+h)-f(x)\|&\geq\|f(y)-f(x)\|-\|f(x+h)-f(y)\|\\
&\geq (md(f_{n,j},x)-\e)(y-x)-\e h\\
&\geq ((md(f_{n,j},x)-\e)(1-\e h\cdot(2n)^{-1}-\e)h.
\end{split}
\end{equation*}
Thus $md_+(f,x)=md(f_{n,j},x)=mD^+(f,x)$. 
\par
A similar argument shows that $md_-(f,x)=md(f_{n,j},x)=mD^+(f,x)$
for $x\in E_{n,j}$,
and thus $md(f,x)$ exists for all $x\in A\setminus N$.
\end{proof}

Theorem~\ref{butterflythm} has the following corollary.

\begin{corollary}\label{butterflycor}
Let $f:\R\to X$ be arbitrary. Then there
exists a set $N\subset\R$ with $\lambda(N)=0$, such that 
if $x\in\R\setminus N$,
and $\min(mD^-(f,x),mD^+(f,x))<\infty$, then $md(f,x)$ exists.
\end{corollary}

Corollary~\ref{butterflycor} together with~\cite[Theorem~2.6]{Dkir}
imply the following:

\begin{corollary}\label{mdifbutterflycor}
Let $f:\R\to X$ be arbitrary. Then there
exists a set $M\subset\R$ with $\lambda(M)=0$, such that if 
$x\in\R\setminus M$,
and $\min(mD^-(f,x),mD^+(f,x))<\infty$, then $f$ is metrically differentiable
at~$x$.
\end{corollary}

\section{Points of metric non-differentiability}\label{ptsmetrnondifsec}

We will use following lemma proved in~\cite[Lemma~2.4]{DZ2}.

\begin{lemma}\label{mdflem} Let $f:[c,d]\to X$, $x\in[c,d]$. Then the following hold.
\begin{enumerate}
\item If $md(f,x)=0$, then $f$ is metrically differentiable at $x$.
\item If $h:[a,b]\to[c,d]$ is differentiable at $w\in[a,b]$,
$h(w)=x$, and $f$ is metrically differentiable at~$x$,
then
$f\circ h$ is metrically differentiable at~$w$,
and $md(f\circ h,w)=md(f,x)\cdot|h'(w)|$.
\end{enumerate}
\end{lemma}

\begin{lemma}\label{kirfedlem} 
Let $X$ be a Banach space, and let $f:[a,b]\to X$. If $md(f,\cdot)$
is continuous at $x\in[a,b]$,
then there exists $\delta>0$ such that 
\[ \bigvee^t_s f=\int^t_s md(f,y)\,dy\quad\text{ for all }s<t,s,t\in[x-\delta,x+\delta]\cap[a,b].\]
\end{lemma}

\begin{proof} Let $\delta>0$ be chosen such that for all 
$s\in[x-\delta,x+\delta]\cap[a,b]$ we have that 
$md(f,s)$ exists and
$|md(f,x)-md(f,s)|\leq1$.
It follows from~\cite[\S 2.2.7]{F} that $f|_{[x-\delta,x+\delta]\cap[a,b]}$ is Lipschitz.
We obtain that
\[ \int^t_s md(f,y)\,dy=\int_{f([s,t])} N(f|_{[s,t]},y)\,d\Hau^1(y)=\bigvee^t_s f,\]
for all $s<t$, $s,t\in[x-\delta,x+\delta]\cap[a,b]$
(here, $N(f|_{[s,t]},y)$ is the multiplicity with which
the function $f|_{[s,t]}$ assumes a value $y$). The first equality
follows from~\cite[Theorem~7]{K}, the second equality follows from~\cite[Theorem~2.10.13]{F}.
\end{proof}

Let $f:[a,b]\to X$, $I=[a,b]$.
We say that $x\in I$ is {\em metrically regular point of the function~$f$}, provided
\[ \lim_{\substack{t\to0\\x+t\in I}} 
\frac{\|f(x+t)-f(x)\|}{\big|\bigvee^{x+t}_{x} f\big|}=1.\]

\begin{lemma}\label{ctstomreglem} 
Let $X$ be a Banach space, $g:[a,b]\to X$, $x\in[a,b]$, 
$md(g,x)>0$, and $md(g,\cdot)$ is continuous at $x$. Then $x$ is metrically 
regular point of the function~$g$.
\end{lemma}

\begin{proof}
Let $\e>0$. Find 
$\delta_0>0$ such that 
$(1-\e)\,md(g,x)|t|\leq\|g(x+t)-g(x)\|$
and 
$md(g,x+t)<(1+\e)\cdot md(g,x)$,
whenever $|t|<\delta_0$ and $x+t\in[a,b]$.
Using Lemma~\ref{kirfedlem}, we can find $0<\delta<\delta_0$
such that for all $|t|<\delta$ we have 
\begin{equation*}
\begin{split}
\bigg(\frac{1-\e}{1+\e}\bigg)\bigg| \bigvee^{x+t}_{x}g\bigg| &=
\bigg(\frac{1-\e}{1+\e}\bigg) \bigg|\int^{x+t}_x md(g,s)\,ds \bigg|\\
&\leq(1-\e)\cdot md(g,x)\,|t|\leq\|g(x+t)-g(x)\|
\leq\bigg|\bigvee^{x+t}_{x}g\bigg|.
\end{split}
\end{equation*}
If $t\neq0$, by dividing by $\big|\bigvee^{x+t}_{x}g\big|$ (which is strictly positive),
we obtain
$\frac{1-\e}{1+\e}\leq\frac{\|g(x+t)-g(x)\|}{|\bigvee^{x+t}_{x}g|}\leq 1$,
and thus $x$ is metrically regular point of $f$.
\end{proof}

The following lemma shows that the condition~\eqref{md} is satisfied
``unilaterally'' at a point~$x$ provided $md(f,\cdot)$ is continuous
at~$x$.

\begin{lemma}\label{ctstometric}
Let $X$ be a Banach space, and let $f:[a,b]\to X$. If $md(f,\cdot)$
is continuous at $x\in[a,b]$, then 
\begin{equation}\label{onesidedmetrdif} 
\|f(y)-f(z)\|-md(f,x)|y-z|=o(|x-z|+|x-y|),
\end{equation}
whenever $(y,z)\to(x,x)$, and $\sign(z-x)=\sign(y-x)$.
\end{lemma}

\begin{proof}
If $md(f,x)=0$, then the conclusion follows from Lemma~\ref{mdflem},
so we can assume that $md(f,x)>0$.
Lemma~\ref{ctstomreglem} implies that $x$ is metrically regular
point of $f$.
Now we will prove that $f$ satisfies~\eqref{onesidedmetrdif} at $x$.
Let $0<\e<1$.
Using Lemma~\ref{kirfedlem}, find $\delta>0$ such that for all $t$ with
$x+t\in[a,b]\cap[x-\delta,x+\delta]$
we have that
$(1-\e)\,\big|\bigvee^{x+t}_x f\big|\leq \|f(x+t)-f(x)\|$, 
\[(1-\e)\,md(f,x)\leq md(f,x+t)\leq (1+\e)\,md(f,x),\]
and
$\bigvee^z_y f=\int^z_y md(f,s)\,ds$
for all $y,z\in[a,b]\cap[x-\delta,x+\delta]$.
Let $y,z\in[a,b]\cap[x-\delta,x+\delta]$ with $\sign(z-x)=\sign(y-x)$.
Without any loss of generality, we can assume
that $z>x$, and $|z-x|\geq|y-z|$.
We obtain that 
\begin{equation*}
\begin{split}
\|f(y)-f(z)\| &\geq \|f(z)-f(x)\|-\|f(y)-f(x)\|\\
& \geq (1-\e)\bigvee^z_x f-\|f(y)-f(x)\|\\
& \geq (1-\e)\int^z_x md(f,t)\,dt-\bigvee^y_x f\\
& \geq (1-\e)\int^z_x md(f,t)\,dt-\bigg|\int^y_x md(f,t)\,dt\bigg|\\
& \geq (1-\e)^2\,md(f,x)(z-x)-(1+\e)\,md(f,x)|y-x|\\
&= md(f,x)|z-y|-\e\cdot\underbrace{((2-\e)(z-x)+|y-x|)\cdot md(f,x)}_{\eta(\e,y,z)}.
\end{split}
\end{equation*}
It is easy to see that $\frac{\eta(\e,y,z)}{|z-x|+|y-x|}$ is bounded 
from above by $2\cdot md(f,x)$ for all
$\e\in(0,1)$. 
For the other inequality, note that
\begin{equation}\label{eqreus} \|f(y)-f(z)\|\leq\bigvee^z_y f=\int^z_y md(f,s)\,ds\leq
(1+\e)\,md(f,x)|z-y|,\end{equation}
and the conclusion easily follows.
\end{proof}

We now show that if the metric derivative of $f$ exists at each point
and is continuous, 
then the mapping is metrically differentiable on a large set of
points. We prove this in several steps.

\begin{proposition}\label{nonmdprop}
Let $X$ be a Banach space, $f:[a,b]\to X$
be such that $md(f,x)=1$ for each $x\in[a,b]$.
Then the set of points $x\in[a,b]$ such that 
$f$ is not metrically differentiable at $x$,
is $\sigma$-symmetrically porous.
\end{proposition}

\begin{proof} 
Let $A$ be the set of points $x\in(a,b)$ such that 
$f$ is not metrically differentiable at $x$.
By Lemma~\ref{ctstometric}, we see that the condition~\eqref{md}
is satisfied unilaterally at each $x\in[a,b]$.
\par
Suppose that $x\in A$.
We claim that there exist $\delta_j=\delta_j(x)\to0+$ such that
\begin{equation}\label{symm}
\liminf_{j\to\infty} \frac{\|f(x+\delta_j)-f(x-\delta_j)\|}{2\delta_j}<1.
\end{equation}
To see this, note that because $x\in A$, there exist $(y_j)_j$, $(z_j)_j$ such that
$y_j<x<z_j$ (because~\eqref{md} is satisfied unilaterally at~$x$), $\lim_j y_j=\lim_j z_j=x$, and
$
\liminf_{j\to\infty} \frac{\|f(z_j)-f(y_j)\|}{z_j-y_j}<1-\e,
$
for some $\e>0$.
Without any loss of generality, we can assume that $z_j-x\leq x-y_j$.
Let $\tilde y_j=2x-y_j$, and note that $z_j\leq\tilde y_j$. 
If $\tilde y_j=z_j$, take $\delta_j=z_j-x$, otherwise
note that for $j\in\N$ large enough we have
\begin{equation*}
\begin{split}
\|f(y_j)-f(\tilde y_j)\|&\leq \|f(y_j)-f(z_j)\|+\|f(z_j)-f(\tilde y_j)\|\\
&\leq(1-\e)(z_j-y_j)+(\tilde y_j-z_j).
\end{split}
\end{equation*}
Now, as $\tilde y_j-z_j\leq z_j-y_j$, we obtain
$(\tilde y_j-z_j)-\frac{\e}{2}(z_j-y_j)\leq \big(1-\frac{\e}{2}\big)(\tilde y_j-z_j)$,
and thus
$\|f(y_j)-f(\tilde y_j)\|\leq \big(1-\frac{\e}{2}\big)(\tilde y_j-y_j)$.
Now define $\delta_j=\tilde y_j-x=x-y_j$, and~\eqref{symm} follows.
\par
Let $A_{nm}$ be the set of all $x\in A$ such that
\begin{itemize}
\item
there exist a sequence $(\delta_j)_j$, 
such that $\delta_j\to0+$, 
and 
$\|f(x-\delta_j)-f(x+\delta_j)\|\leq \big(1-\frac{1}{m}\big)\,2\,\delta_j$,
\item 
for each $t\in[0,1]$ with $0<|x-t|<1/n$ we have
$\big(1-\frac{1}{2m}\big)|x-t|<\|f(t)-f(x)\|$.
\end{itemize}
By the above argument,
it is easy to see that $A=\bigcup_{n,m} A_{nm}$.
\par
Fix $n,m\in\N$. Let $x\in A_{nm}$. 
There exists $j_0\in\N$ such
that for all $j\geq j_0$ we have $0<\delta_j<(4n)^{-1}$.
Let $z_j:=x+\delta_j$, and $y_j:=x-\delta_j$.
Fix $j\geq j_0$ and suppose that $w\in [z_j,z_j+2\delta_j]$.
Then $|w-y_j|<1/n$, and we have that 
\begin{equation*}
\begin{split}
\|f(y_j)-f(w)\| & \leq \|f(y_j)-f(z_j)\|+\|f(z_j)-f(w)\|\\
&\leq \bigg(1-\frac{1}{m}\bigg) \,2\delta_j+|w-z_j|.
\end{split}
\end{equation*}
By the choice of $w$ we have 
$-\frac{2\delta_j}{2m}+(w-z_j)\leq\big(1-\frac{1}{2m}\big)(w-z_j)$,
and thus
$\|f(y_j)-f(w)\| \leq \big(1-\frac{1}{2m}\big)(w-y_j)$.
This implies that $w\not\in A_{nm}$. 
We obtained that $[z_j,z_j+2\delta_j]\cap A_{nm}=\emptyset$.
Similarly, $[y_j-2\delta_j,y_j]\cap A_{nm}=\emptyset$,
and the symmetric porosity of $A_{nm}$ follows.
\end{proof}

We will need the following auxiliary lemma.

\begin{lemma}\label{symimlem}
Let $B\subset[a,b]$ be symmetrically porous and $h:[a,b]\to\R$
be continuously differentiable and bilipschitz.
Then $h(B)$ is symmetrically porous.
\end{lemma}

\begin{proof}
Let $L>0$ be such that $L^{-1}|x-y|\leq|h(x)-h(y)|\leq L|x-y|$ 
for all $x,y\in[a,b]$.
Let $x\in B$. 
Let $\delta_n,\alpha_n>0$ be such that $B(x-\delta_n,\alpha_n)\cup
B(x+\delta_n,\alpha_n)\subset\R\setminus B$, $\alpha_n\to 0$, and
$c\delta_n\leq \alpha_n$.
First, we will show that 
\begin{equation}\label{symeq1}
B(h(x\pm\delta_n),\alpha_n/(2L))\cap h(B)=\emptyset.
\end{equation}
Let $z\in B(h(x\pm\delta_n),\alpha_n/(2L))$. If $y\in[a,b]$ is such 
that $h(y)=z$, then 
\[ |x\pm\delta_{{n}}-y|\leq L|h(x\pm\delta_{{n}})-h(y)|
\leq\alpha_n/2,\]
and thus $z\not\in h(B)$ (since $h$ is one-to-one), and~\eqref{symeq1} holds.
\par
Note that since $h'(x)\ne0$,
\[ \left|1-\frac{h(x+\delta_n)-h(x)}{h(x)-h(x-\delta_n)}\right|
=\left|1-\frac{h'(x)\delta_n+o(\delta_n)}{h'(x)\delta_n+o(\delta_n)}\right|\to 
{0},\]
as $n\to\infty$,
and thus 
\begin{equation}
\label{symeq2}
| h(x+\delta_n)-h(x) - (h(x)-h({x-}\delta_n))|
\leq c/(4L^2)|h({x-\delta_n})-h(x)|\leq c\delta_n/(4L)
{\leq}\alpha_n/(4L)
\end{equation}
for $n$ large enough.
Now, we will show that 
$2h(x)-z\in B(h(x+\delta_{n}),\alpha_n/(2L))$ 
whenever $z\in B(h(x-\delta_{{n}}),\alpha_n/(4L))$, and $n$ is large enough.
Together with~\eqref{symeq1},
this easily implies that $S\overline{p}(h(B),h(x))>0$.

Assume that $z\in B(h(x-\delta_{n}),\alpha_n/(4L))$,
and that~\eqref{symeq2} holds.
Then
\begin{equation*}
\begin{split} 
|2h(x)-z-h(x+\delta_{{n}})|&\leq 
|h(x+\delta_{n})-h(x)+(h(x-\delta_{n})-h(x))|
+|h(x-\delta_{n})-z|\\
&\leq \alpha_n/(4L)+\alpha_n/(4L)=\alpha_n/(2L),
\end{split}
\end{equation*}
and thus $2h(x)-z\in B(h(x+\delta_n),\alpha_n/(2L))$, 
and the conclusion follows.
\end{proof}

We have the following:

\begin{theorem}\label{ptmetrnondif}
Let $f:[a,b]\to X$ be such that $md(f,\cdot)$ 
is continuous on $[a,b]$. Then the set of points, where $f$ is not
metrically differentiable, is $\sigma$-symmetrically porous.
\end{theorem}

\begin{proof} Let $A\subset[a,b]$ be the set where $f$ is not
metrically differentiable. Lemma~\ref{mdflem} implies that if $x\in A$,
then $md(f,x)>0$. Let $A=\bigcup_n A_n$, where $A_n=\{x\in A:md(f,x)>1/n\}$.
It is enough to show that each $A_n$ is $\sigma$-symmetrically porous.
Because $md(f,\cdot)$ is continuous, we have that each $A_n$ is open.
Let $(c,d)$ be an open component of $A_n$, let $g=f|_{[c,d]}$,
$G=g\circ v_g^{-1}$
(see Section~\ref{preliminaries} for the definition of $v_g$). 
Using Lemma~\ref{kirfedlem}, it is easy to see that
$md(G,x)=1$ for all $x\in v_g((c,d))$.
Then Proposition~\ref{nonmdprop} implies that 
$G$ is metrically differentiable outside a $\sigma$-symmetrically 
porous set $B$.
Because $v_g$ is continuously differentiable and bilipschitz,
by Lemmas~\ref{mdflem}~and~\ref{symimlem},
we obtain that $g=G\circ v_g$ is metrically differentiable
outside a $\sigma$-symmetrically porous set $v_g^{-1}(B)$.
\end{proof}


\begin{remark}
It is easy to see that if $f$ is a real-valued function and
$md(f,\cdot)$ is continuous on $[a,b]$, then the set of points
where $f$ is not metrically differentiable is at most countable.
However, in Theorem~\ref{difexample} below,
we show that already in a $2$-dimensional situation
such a set may be uncountable.
Thus, Theorem~\ref{difexample} shows that Theorem~\ref{ptmetrnondif} cannot
be strengthened to make the exceptional set countable.
\end{remark}

\begin{theorem}\label{difexample}
For any norm $\|\cdot\|$ in the $2$-dimensional plane,
there exists a curve $\gamma:[0,\ell]\to (\R^2,\|\cdot\|)$ with
$md(\gamma,x)=1$ for all $x\in[0,\ell]$, but such that
the set of points where $\gamma$ is not metrically differentiable
is uncountable.
\end{theorem}

We will give a detailed proof of 
this theorem for $\|\cdot\|$ being the Euclidean norm.
In Remark~\ref{explanation}, we explain how this case reflects the
most general situation. Note however, that if one uses a ``polygonal'' norm 
(for example, the $\ell_1$-norm), then much simpler constructions are
possible. We  explain this in Remark~\ref{explanation-next}.

Before we start the proof of Theorem~\ref{difexample}, let us establish
the following 
property of logarithmic spirals, which will be used in 
the proof of Lemma~\ref{th-logsp}.

\begin{lemma}\label{lem-logsp}
Assume $S_{a,b}$ is a planar curve defined 
in polar coordinates $(r,\phi)$ by the equation $r=ae^{b\phi}$
with $a>0$, $b\ne0$ (logarithmic spiral).
Then
the length of the 
arc of $S_{a,b}$ between the origin
and the point with modulus
$r_0$ and argument $\phi_0$ 
is equal to $\frac{\sqrt{b^2+1}}{|b|}r_0$.

In other words, if $S_{a,b}:[0,+\infty)\to\complex$
is the arc-length parameterization of this logarithmic
spiral such that $S_{a,b}(0)=0$,
then
\begin{equation}
\frac{|S_{a,b}(t)|}{t}=\frac{|b|}{\sqrt{b^2+1}}
\label{eq-logsp}
\end{equation}
for all $t>0$.
\end{lemma}

\begin{proof}
A routine computation of the length of the logarithmic spiral with
the given equation in polar coordinates proves the lemma.
\end{proof}

\begin{lemma}
\label{th-logsp}
For any angle $\alpha\in(0,\pi/2)$ and a constant $\mycon\in(0,1)$
there is a piecewise smooth
planar curve such that its arc-length parameterization
$g=g_{\mycon,\alpha}:\R\to\complex$ has
the following properties:
\begin{enumerate}
\item[(a)] $g([0,1])$ is a horizontal interval 
and there exists $L_{\mycon,\alpha}>0$
such that $g([1+L_{\mycon,\alpha},+\infty))$ 
and $g((-\infty,-L_{\mycon,\alpha}])$ are horizontal rays;
\item[(b)] there exists  $t_{\mycon,\alpha}>1/2$ such that 
the arguments
of $z_{\pm}=g(1/2\pm t_{\mycon,\alpha})-g(1/2)$ are equal to $(-\alpha)$ and
$(\pi+\alpha)$ resp.; 
\item[(c)] $|g(t)-g(s)|/|t-s|>\mycon$
for all  $s\in[0,1]$ and $t\ne s$.
\end{enumerate}
\end{lemma}

\begin{proof}
Let $B>0$ be large enough so as to ensure that 
\begin{align}
\label{conds-B}
\frac{B}{\sqrt{B^2+1}}>\mycon
\quad
\text{ and }
-B \sin\alpha+\cos\alpha<0. 
\end{align}
In~\eqref{conds-B-add-sec},
we will impose another condition on $B$ 
which also bounds $B$ from below.
Fix $b>B$ and
denote $k=\frac{b}{\sqrt{b^2+1}}$.

We first construct a  piecewise smooth planar curve
$f=f_{\mycon,\alpha}:[0,+\infty)\to\complex$ such that 
\begin{equation}
\label{def-g}
g(t)=
\begin{cases}
f(t),&\text{if }t\ge0,\\
1-\overline{f(1-t)},&\text{if }t<0,
\end{cases}
\end{equation}
has the desired properties.

For $t\in[0,1]$ we set $f(t)=t+0i$. 
Now 
let $S_{1,-b}:[0,+\infty)\to\complex$
be the arc-length parameterization of the
logarithmic spiral from Lemma~\ref{lem-logsp}.
Identity \eqref{eq-logsp} implies that
the point $S_{1,-b}(k^{-1})$ 
has modulus $1$, therefore, it coincides with
$f(1)$.

For $t\in[1,1+k^{-1}(e^{b\alpha}-1)]$ 
we put
$f(t)=S_{1,-b}(t+k^{-1}-1)$.
Then for every $s\ge0$ one has:
\begin{equation}
\label{prop-f-1}
\frac{1+s}{|f(1+s)-f(0)|}
<
\frac{k^{-1}+s}{|f(1+s)-f(0)|}
=
\frac{k^{-1}+s}{|S_{1,-b}(k^{-1}+s)|}
=k^{-1}.
\end{equation}

Let $s_0=k^{-1}(e^{b\alpha}-1)$. 
Then the point $f(1+s_0)=S_{1,-b}(k^{-1}e^{b\alpha})$ 
has  modulus $e^{b\alpha}$ and argument $-\alpha$.

Now let
$S_{e^{2b\alpha},b}:[0,+\infty)\to\complex$
be another logarithmic spiral parametrized by the arc-length.
For $t\in[1+s_0,1+s_0+s_1]$ 
(where $s_1$ is defined below),
let $f(t)=S_{e^{2b\alpha},b}(t+k^{-1}-1)$. 
Again, note that 
$S_{1,-b}(t+k^{-1}-1)$ and
$S_{e^{2b\alpha},b}(t+k^{-1}-1)$
are equal at $t=1+s_0$, since by \eqref{eq-logsp}
the lengths of the arcs of both logarithmic spirals
between the origin and the point 
with modulus $e^{b\alpha}$ and argument $-\alpha$
are equal to $k^{-1}e^{b\alpha}=k^{-1}+s_0$.
Furthermore, for every $s\ge0$ one has:
\begin{equation}
\label{prop-f-2}
\frac{1+s_0+s}{|f(1+s_0+s)-f(0)|}
<
\frac{k^{-1}+s_0+s}{|S_{e^{2b\alpha,},b}(k^{-1}+s_0+s)|}
=k^{-1}.
\end{equation}

Let us find the slope of the tangent to the logarithmic spiral
$S_{e^{2b\alpha},b}$ at the point with modulus
$e^{b\alpha}$ and argument $-\alpha$. If we denote by
$z(\phi)=e^{2b\alpha}e^{b\phi}e^{i\phi}$ the polar parameterization
of $S_{e^{2b\alpha},b}$,
then
$\text{Im}\,  \frac{dz}{d\phi}(-\alpha)$
is equal to
$e^{b\alpha}(-b \sin\alpha+\cos\alpha)<0$. 
Therefore, the $y$-coordinate of $f(t)$ continues to decrease
as $\phi$ increases from $-\alpha$ to some $-\beta\in(-\alpha,0)$ 
such that $-b \sin\beta+\cos\beta=0$
(i.e., $\tan\beta=1/b$). 
Let $s_1$ be such that
$f(1+s_0+s_1)=S_{e^{2b\alpha},b}(e^{2b\alpha}e^{-b\beta})$ 
is the point
with modulus $e^{2b\alpha-b\beta}$ and argument $-\beta$,
i.e.,
$s_1=k^{-1}e^{b\alpha}(e^{b(\alpha-\beta)}-1)$. 

For $t\ge1+s_0+s_1$, define $f(t)$ as $f(1+s_0+s_1)+(t-1-s_0-s_1)$. Then
one easily checks that since $\cos\beta=k$, the law of
cosines for the triangle with vertices
in $f(0)$, $f(1+s_0+s_1)$ and $f(1+s_0+s_1+s)$ guarantees that the 
inequality
\begin{equation}
\label{prop-f-3}
\frac{(1+s_0+s_1+s)^2}{|f(1+s_0+s_1+s)-f(0)|^2}
\le
k^{-2}
\end{equation}
holds for all $s\ge0$.

Inequalities \eqref{prop-f-1}, \eqref{prop-f-2}, \eqref{prop-f-3}
imply that 
\begin{equation}
\label{prop-f}
\frac{t}{|f(t)-f(0)|}\le k^{-1}
\end{equation}
for all $t>0$.

Note that if we now define $g$ as in \eqref{def-g}, then property $(a)$ in
the Lemma holds for $L_{\mycon,\alpha}=s_0+s_1$.

The argument of
$f(1+s_0)-f(0)$ is equal to $(-\alpha)$. 
Then
the arguments
of $g(1/2\pm (1/2+s_0))-g(1/2)$ are equal to $(-\alpha')$ and
$(\pi+\alpha')$ respectively, where $\alpha'>\alpha$.
Since the argument of $g(1/2+ t)-g(1/2)$ is continuous
in $t$, there is a value $t_{\mycon,\alpha}$ between $1/2$ and
$1/2+s_0$ such that property $(b)$ in the present lemma holds for
$t_{\mycon,\alpha}$.

We have already proved, see \eqref{prop-f}, that
property $(c)$ in the present lemma holds
for $s=0$ and all $t>0$ (as $g(t)=f(t)$ for $t\ge0$). 
If $s\in[0,1]$ and $1\le t\le 1+s_0+s_1$, then $g(s)=s$ and
\begin{align*}
t-s
&=(t-1)+(1-s)
=(k^{-1}|g(t)|-k^{-1})+(1-s)\\
&=k^{-1}(|g(t)|-s)-(k^{-1}-1)(1-s)
<k^{-1}|g(t)-g(s)|.
\end{align*}
If $t\ge 1+s_0+s_1$, then
$|g(t)-g(0)|/t\ge k$, and therefore,
$$
\frac{|g(t)-g(s)|}{t-s}
\ge
\frac{|g(t)|-s}{t-s}
\ge
\frac{k-x}{1-x},
$$
where $x=\frac{s}{t}\le\frac{1}{1+s_0+s_1}$.
Then $\frac{k-x}{1-x}\ge k-\frac{1-k}{s_0+s_1}
\ge k-\frac{k(1-k)}{e^{b\alpha}-1}$.
Note that the latter expression is an
increasing function of $b$ (as $k$ is a function 
of $b$), which tends to $1$ as $b$ tends to infinity. Therefore, if 
in addition to \eqref{conds-B} 
we require that
\begin{equation}
\label{conds-B-add-sec}
\frac{B}{\sqrt{B^2+1}}-\frac{\frac{B}{\sqrt{B^2+1}}(1-\frac{B}{\sqrt{B^2+1}})}
{e^{B\alpha}-1}
>\mycon,
\end{equation}
then property $(c)$ in the Lemma holds for all $s\in[0,1]$ and $t\ge1$. 
It remains to note that this property trivially holds for $s,t\in[0,1]$
and that by symmetry, 
the case $s\in[0,1]$, $t<0$ is analogous to $1-s\in[0,1]$,
$1-t>1$.

Thus, conditions $(a)$--$(c)$ hold for $g$ with 
$t_{\mycon,\alpha}\in(1/2,1/2+s_0)$
and $L_{\mycon,\alpha}=s_0+s_1$.
\end{proof}
\begin{remark}
\label{rem-Lip}
In addition to properties $(a)$--$(c)$
of Lemma~\ref{th-logsp} we may 
assume that 
the curve $g_{\mycon,\alpha}$
is a 
graph of Lipschitz piecewise smooth
function $F_{\mycon,\alpha}:\R\to\R$.
\end{remark}
\begin{proof}
Let us analyze the tangent vector to $g_{\mycon,\alpha}$
when the argument 
$t$ changes from $1$ to $1+s_0$ and from $1+s_0$ to $1+s_0+s_1$
(see the proof of Lemma~\ref{th-logsp}).

The arc of $g_{\mycon,\alpha}$ between $g_{\mycon,\alpha}(1)$ 
and $g_{\mycon,\alpha}(1+s_0)$ 
has the polar parameterization
$z(\phi)=e^{b\phi}e^{-i\phi}$, $\phi$ increases from $0$ to $\alpha$.
Then the $x$-coordinate 
$\text{Re}\,  \frac{dz}{d\phi}(\phi)$
of the tangent vector
is equal to
$e^{b\phi}(b \cos\phi-\sin\phi)$.
This is positive provided $\tan\phi<b$. Thus, 
we impose the following additional restriction on $B$:
\begin{equation}
\tan \alpha < B.
\end{equation}
Since 
$\text{Re}\,  \frac{dz}{d\phi}(\phi)$
is continuous, we conclude that its minimum on 
$\phi\in[0,\alpha]$ is positive.
The $y$-coordinate of the tangent vector is continuous in $t$,
therefore is bounded for $t\in[1,1+s_0]$.
Thus, the slope of the tangent vector is bounded.
Hence $g_{\mycon,\alpha}|_{[1,1+s_0]}$ is a graph of Lipschitz
function.

The arc of $g_{\mycon,\alpha}$ between $g_{\mycon,\alpha}(1+s_0)$ 
and $g_{\mycon,\alpha}(1+s_0+s_1)$ 
has the polar parameterization
$z(\phi)=e^{2b\alpha+b\phi}e^{i\phi}$, $\phi$ increases from $-\alpha$
to $-\beta$.
Then
$\text{Re}\,  \frac{dz}{d\phi}(\phi)
=e^{2b\alpha+b\phi}(b \cos\phi-\sin\phi)>0$
since $\cos\phi>0$ and $\sin\phi<0$ for $\phi\in(-\alpha,-\beta)$.
In the same way this implies that 
 $g_{\mycon,\alpha}|_{[1+s_0,1+s_0+s_1]}$ is a graph of Lipschitz
function.
\end{proof}

\begin{proof}[Proof of Theorem~\ref{difexample} in the Euclidean case]
Let $\alpha_n\to\pi/2$ and $\mycon_n\to1$, 
$n\ge1$ be two increasing sequences
of positive reals.
For every pair $(\alpha_n,\mycon_n)$
consider a Lipschitz function $F_n(x)=F_{\mycon_n,\alpha_n}(x+1/2)$, where
$F_{\mycon_n,\alpha_n}$
is a Lipschitz piecewise smooth function
described in Remark~\ref{rem-Lip} (whose graph is the
curve $g_{\mycon_n,\alpha_n}$ from Lemma~\ref{th-logsp}).
The function $F_n$ is even.
Note that $F_n(x)$ is constant for
$|x|\ge x_n=\text{Re}\,(g_{\mycon_n,\alpha_n}(1+L_{\mycon_n,\alpha_n}))$.
Denote by $\Gamma_{F_n}(x)=x+iF_n(x)$ the graph of
$F_n$ and for each $n\ge1$ choose $L_n>x_n$ such that
\begin{equation}
\label{prop-Ln}
L_n-x_n>n\mathcal H^1 (\Gamma_{F_n}[-x_n,x_n]).
\end{equation}

Now let $G_n(x)=\frac{F_n(L_nx)-F_n(L_n)}{L_n}$.
The function $G_n$ has the following properties:
\begin{itemize}
\item $G_n$ is a nonnegative even piecewise smooth 
Lipschitz function on $\R$, 
\item $G_n$ is zero on $(-\infty,-1]\cup[1,\infty)$,
\item $G_n(x)$ attains its maximum at $x=0$, $G_n(0)<1$ and
$G_n^{-1}(G_n(0))=[-1/L_n,1/L_n]$
($L_n>|F_n(x_n)|$ from \eqref{prop-Ln}),
\item If $\gamma_n=\Gamma_{G_n}$ is the graph of $G_n$, then
$\mathcal H^1( \gamma_n[-a_n,a_n])<1/n$, where 
$a_n=\sup\{x\colon G_n(x)>0\}$,
\item There exists $t_n\in(0,a_n)$ such that 
the argument of 
$\gamma_n(t_n)-\gamma_n(0)$
is equal to $(-\alpha_n)$,
\item 
The ratio 
$\frac{|\gamma_n(x)-\gamma_n(y)|}{\mathcal H^1(\gamma_n[x,y])}$
is bounded from below by $\mycon_n$
for all pairs of $x\ne y$ such that $|x|\le1/L_n$.
\end{itemize}

Denote by $p_n$ the length of $\gamma_n[-a_n,a_n]$.
Let $\gap_n\searrow0$ ($n\ge1$) be such that 
\[
\gap_{n+1}p_{n+1}<(\gap_np_n)/4, \quad
2\gap_{n+1}<\gap_n/L_n, \quad
\gap_1<1/2, \text{ and }
\sum_{n\ge1}\gap_n<1.
\]
The first property of $\gap_n$ guarantees that for every $n\ge1$
\[
\sum_{j\ge1}2^{j-1}\gap_{n+j}p_{n+j}
<2\gap_{n+1}p_{n+1}.
\]
Note that as $G_n(x)$ is a hat-like function on $[-1,1]$, 
the graph of  $G_n^{(\dd,\gap)}(x)=
\gap G_n(\gap^{-1}(x-\dd))$ 
is the rescaled ``hat'' on $[\dd-\gap,\dd+\gap]$.
For any closed interval $I=[\dd-\gap,\dd+\gap]$ denote by
$I^{(L)}$ the interval $[\dd-\gap/L,\dd+\gap/L]$.

\begin{figure}%
\setlength{\unitlength}{\textwidth}%
\begin{picture}(1,0.15)
\put(0.05,0){\resizebox{0.9\textwidth}{!}{%
\includegraphics{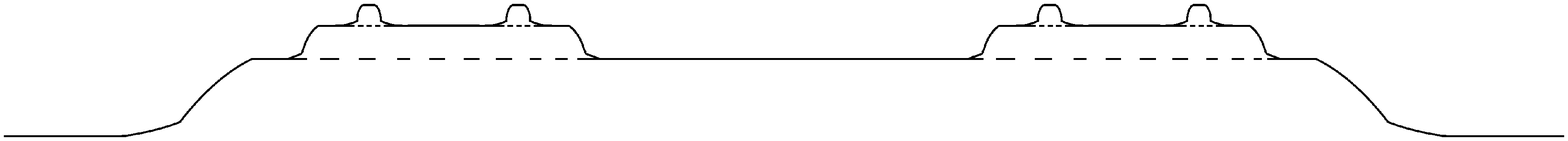}}}%
\put(0.057,0.0149){\circle*{0.005}}%
\put(0.938,0.0149){\circle*{0.005}}
\put(0.03,-0.03){$-1$}%
\put(0.9,-0.03){$-1+2\gap_1$}%
\end{picture}%
\caption{A graph of $S_3(x)$}%
\label{fig}%
\end{figure}

For $x\in[-1,1]$, let
\[
\h_1(x)
=S_1(x)
=\sum_{\dd\in\{-1+\gap_1,1-\gap_1\}}G_1^{(\dd,\gap_1)}(x).
\]
Let 
$\mathcal G_1=\{[-1,-1+2\gap_1],[1-2\gap_1,1]\}$
(since $\gap_1<1/2$, these intervals are disjoint).
Now we define inductively two sequences of families  
of intervals as follows:
\begin{align*}
\mathcal F_n
&=\{I^{(L_n)}
\text{ such that }
I\in\mathcal G_n\};\\
\mathcal G_{n+1}
&=\{[a,a+2\gap_{n+1}],[b-2\gap_{n+1},b]
\text{ such that }
[a,b]\in\mathcal F_n\}.
\end{align*}

For every $n\ge1$, $x\in[-1,1]$, let
\begin{align}
\label{def-Sn}
\h_{n+1}(x)
&=
\sum_{[a,b]\in\mathcal F_{n}}\quad
\sum_{\dd\in\{a+\gap_{n+1},b-\gap_{n+1}\}}
G_{n+1}^{(\dd,\gap_{n+1})}(x);\\
\notag
S_{n+1}(x)&=S_{n}(x)+\h_{n+1}(x)
\end{align}
(Figure~\ref{fig} shows a possible graph of $S_3(x)$).
Note that the definition of $\h_1$ agrees with \eqref{def-Sn} if
we let $\mathcal F_0=\{[-1,1]\}$.
For all $n$, $\mathcal F_n$ consists of
$2^n$ disjoint closed intervals of
the same length $2\gap_n/L_n$, whose union
is equal to the preimage $S_n^{-1}(\max_x S_n(x))$.
Since $4\gap_{n+1}\le2\gap_n/L_n$, intervals in $\mathcal G_{n+1}$
are disjoint.

For $x\in[-1,1]$, define $G(x)=\lim_{n} S_n(x)$.
Note that each $S_n$ is continuous and $|G-S_n|=|\sum_{k\ge n+1} \h_k| 
\le\sum_{k\ge n+1} \gap_k$ which tends to zero as $n\to\infty$.
Therefore, $G$ is continuous. Since the length $\ell$ of the graph 
of $G$ is finite
(it is bounded from above by 
$1+\sum_{n\ge1}2^n\gap_n p_n<1+4\gap_1 p_1<5$),
we conclude that the graph of $G$ has an arc-length parameterization.

Let 
$\gamma=\Gamma_G\colon[-1,1]\to\complex$
be the graph of $G$. The curve
$\gamma$ consists of points of two types:
points in $A_1=\bigcup \Gamma_{S_n}[-1,1]\cap \gamma[-1,1]$ and points in  
$A_2=\gamma[-1,1]\setminus A_1$.
The set $A_2$ is a Cantor-like set which will be described below.

For any $t\in\gamma^{-1}(A_1)$,
the metric  derivative of the normal parameterization of $\gamma$ at $t$
is clearly equal to $1$,
since the functions $S_n$ are piecewise smooth.
Consider $c\in C=\gamma^{-1}(A_2)$. Since $\gamma(c)$ does not belong
to $\Gamma_{S_n}[-1,1]$ for any $n$, there is a sequence of intervals
$I_n\in\mathcal G_n$ such that $c=\bigcap_{n\ge1} I_n$.
Then $\gamma(c)$ corresponds
to a certain infinite sequence $\e\in\{0,1\}^\infty$: 
depending whether $I_{n}$ has center at
$\dd=a+\gap_n$ or at $\dd=b-\gap_n$ (see \eqref{def-Sn}),
we let $\e_{n}$ be equal to $0$ or $1$.
Therefore,
$C$ is a Cantor set, and thus it is uncountable. 
We show that for any $c\in C$, the metric derivative
of the normal parameterization of $\gamma$ 
at $c$ is equal to $1$, but the normal parameterization of $\gamma$ 
is not metrically differentiable at $c$.

For any point $c\in C$
there is a pair of sequences of points
$y_n,z_n\to c$, $y_n<c<z_n$ such that
$G(y_n)=S_n(y_n)$, $G(z_n)=S_n(z_n)$ and
the points $\gamma(y_n)$, $\gamma(z_n)$ and 
$\Gamma_{S_n}(\frac{y_n+z_n}{2})$ form
an isosceles triangle with vertex angle $\pi-2\alpha_n$.
This means that not only  the ratio 
between the distance $|\gamma(y_n)-\gamma(z_n)|$ divided by
the length of $\gamma[y_n,z_n]$
does not tend to $1$, but moreover, it tends to $0$.
Therefore, the normal parameterization of
$\gamma$ is not metrically differentiable at 
$c$.

It remains to show that for any point $c\in C$
the metric derivative of the normal parameterization of
$\gamma$ at $c$ is equal to $1$.
We will show that the ratio
\begin{equation}
\label{ratio}
\mathcal H^1(\gamma[c,c+t])/|\gamma(c+t)-\gamma(c)|
\end{equation}
tends to $1$ as $t\to0$.

Assume $t>0$ is small.
Let $\e\in\{0,1\}^{\infty}$ be a sequence corresponding to $\gamma(c)$. 
Without any loss of generality we may assume 
$c+t\in \bigcup_{I\in\mathcal G_1} I$.
Let $\delta\in\{0,1\}^\infty$ be a sequence corresponding 
to $\gamma(c+t)$. If $\gamma(c+t)\in A_1$, then $\delta$ is a 
finite sequence; otherwise, $\delta$ is infinite.

Since $t$ is small, we may assume that 
$\delta_1=\e_1$. 
Let $n\ge1$ be such that $(\e_1,\dots,\e_n)=
(\delta_1,\dots,\delta_n)$ and $\e_{n+1}\ne\delta_{n+1}$
(if such $n$ does not exist, that is, if the sequence
$\delta$ constitutes the beginning of the infinite 
sequence $\e$, we let $n$ be equal to the length of $\delta$).
Note that when $c$ is fixed and $t$ tends to $0$,
then $n$ tends to $\infty$.

In order to find an upper bound for \eqref{ratio}, we will 
use the following estimate:
\begin{align}
\label{estimate}
\frac{\mathcal H^1(\gamma[c,c+t])}{|\gamma(c+t)-\gamma(c)|}
&\le
\frac{\mathcal H^1(\gamma[x,c+t])+\mathcal H^1(\gamma[c,x])}
{|\gamma(c+t)-\gamma(x)|-|\gamma(x)-\gamma(c)|}\\
\notag
&\le
\Bigl(\frac{\mathcal H^1(\gamma[x,c+t])}{|\gamma(c+t)-\gamma(x)|}+y\Bigr)/
(1-y),
\end{align}
for any $x\in(c,c+t)$,
such that the expression $y=\frac{\mathcal H^1(\gamma[c,x])}{|\gamma(c+t)-\gamma(x)|}$
is strictly less than $1$.

Consider first the case when $\delta$ coincides with 
$(\e_1,\dots,\e_n)$. 
In this case,
$G(c+t)=S_{n}(c+t)$ and 
there is an interval 
$I_n\in\mathcal G_{n}$ of length $2\gap_n$ 
containing both $c$ and $c+t$.
Let $J_{1},J_{2}\subset I_n$ be disjoint intervals in
$\mathcal G_{n+1}$ such that $c\in J_{1}\cup J_2$. Since $\delta$
has length $n$, we get $c+t\not\in J_{1}\cup J_{2}$.
Also note that 
$c\in J_1^{(L_{n+1})}\cup J_2^{(L_{n+1})}$
since $G(c)\ne S_{n+1}(c)$.

If $c\in J_i^{(L_{n+1})}$, then let 
$x=\sup\{z\in J_i\colon S_{n+1}(z)>S_n(z)\}$.
Since $S_{n-1}|_{I_n}$ is constant and $G(x)=S_n(x)$, $G(c+t)=S_n(c+t)$,
we may deduce that 
by the property of $G_{n+1}$, the expression
$\frac{\mathcal H^1(\gamma[x,c+t])}{|\gamma(c+t)-\gamma(x)|}$
does not exceed $\mycon_{n+1}^{-1}$.
Now we want to find an upper estimate for
$y=\frac{\mathcal H^1(\gamma[c,x])}{|\gamma(c+t)-\gamma(x)|}$.
The numerator is not greater than
$\sum_{j\ge1}2^{j-1}\gap_{n+j}p_{n+j}<2\gap_{n+1}p_{n+1}$, and the
denominator is at least $\gap_{n+1}(n+1)p_{n+1}$ (this follows from the 
property of $L_n$, see \eqref{prop-Ln}). Therefore, $y\le 2/(n+1)$.
Thus, the quantity \eqref{ratio} is at most 
$\psi_{n+1}(\mycon_{n+1}^{-1})$, where
$\psi_k(t)=(t+2/k)/(1-2/k)$.

Now consider the case when $\delta$ has length at least $n+1$. 
In the above notation this implies that $c\in J_1^{(L_{n+1})}$
and $c+t\in J_2$. Choose
$x=\sup\{z\in J_1\colon S_{n+1}(z)>S_n(z)\}$ as before.
If $c+t\in J_2^{(L_{n+1})}$, then the same proof as in the previous paragraph
shows that 
$\frac{\mathcal H^1(\gamma[x,c+t])}{|\gamma(c+t)-\gamma(x)|}\le
\psi_{n+1}(1)$ (in this case $\gamma$ connects $\gamma(x)$ and 
$\gamma(x')$, where $x'=\inf\{z\in J_2\colon S_{n+1}(z)>S_n(z)\}$, by
a straight line interval).
Thus, the quantity \eqref{ratio} is at most 
$\psi_{n+1}(\psi_{n+1}(1))$.

If $c+t\in J_2\setminus J_2^{(L_{n+1})}$, then 
$\mathcal H^1(\gamma[x,c+t])\le|\gamma(c+t)-\gamma(x)|+\sum_{j\ge1}
2^{j-1}\gap_{n+j}p_{n+j}<|\gamma(c+t)-\gamma(x)|+2\gap_{n+1}p_{n+1}$,
so together with $|\gamma(c+t)-\gamma(x)|>\gap_{n+1}(n+1)p_{n+1}$
we get that the quantity \eqref{ratio} is at most $\psi_{n+1}(1+2/(n+1))$.

It remains to observe that 
the length $n$ of the initial part of sequences $\e$ and $\delta$ tends to
$\infty$ as $t\to0$ and to note that
$\psi_{n+1}(\mycon_{n+1}^{-1})$,
$\psi_{n+1}(\psi_{n+1}(1))$ and
$\psi_{n+1}(1+2/(n+1))$ tend to $1$ as $n$ tends to infinity.
\end{proof}

\begin{remark}
Note that in fact we proved that the curve $\gamma$ constructed above has 
the following property: for every $c\in C$ there exist 
$y_n<c<z_n$ such that $(y_n,z_n)\to(c,c)$ and
\[\frac{|\gamma(y_n)-\gamma(z_n)|}{\mathcal H^1(\gamma[y_n,z_n])}
\to0.\]
This means that this curve has uncountably many ``spikes''.
\end{remark}

\begin{remark}\label{explanation}
For a general norm $\|\cdot\|$ on the $2$-dimensional plane, one can
produce an analogue of the curve constructed in Lemma~\ref{th-logsp}
in the following way.

We may assume the $\|\cdot\|$-norm of the point $1$ on the
complex plane is equal to $1$.
Define $g([0,1])$ to be a horizontal interval as in
\eqref{def-g} ($f(t)=t$ for $t\in[0,1]$), then $\|g(t)\|=t$
for $0\le t\le1$.
Next find a small $\e>0$, such that if we define
$g(t)|_{t>1}$ to be a ray with slope 
$-\e$, then
the condition $(c)$ in Lemma~\ref{th-logsp} with the norm $\|\cdot\|$
instead of Euclidean norm $|\cdot|$ holds for all $t>1$.
Next thing would be to note that the ratio $\|g(t)-g(s)\|/|t-s|$ tends to 
$1$ as $s$ remains in $[0,1]$ and $t$ tends to infinity
(we define $g(t)=1+(t-1)z_{-\e}$, where $\|z_{-\e}\|=1$
and $\tan\arg z_{-\e}=-\e$). So we may choose
a sufficiently large $T_1$ such that if we redefine
$g(t)|_{t>T_1}$ to be a ray with slope $-2\e$, 
then we again have 
condition $(c)$ in Lemma~\ref{th-logsp} 
still valid for
$\|\cdot\|$.
If we continue this way, the curve $g$ would consist of straight intervals
such that each new interval ``turns'' by less than $-\e$ 
with respect to the previous 
interval, and in the end point of each interval the ratio from
condition $(c)$ is very close to $1$ (much closer to $1$ than $\mycon$ is).
Since $N\e\to\infty$, the angle between the horizontal axis and the
subsequent intervals which form the curve $g$ tends to $\pi/2$. So 
there will be a moment when this angle 
becomes bigger than $\alpha$.
At this moment, we stop the process, and start ``rotating'' intervals 
towards horizontal axis (making the slope less negative)
in order to obtain a broken line
satisfying the conditions $(a)$--$(c)$ from Lemma~\ref{th-logsp}.

One can check that since the arc-length parameterization of the
boundary of a unit ball of arbitrary norm
is uniformly continuous, the 
algorithm explained above can be implemented for every $2$-dimensional norm
(of course, $\e$ would depend on the norm).

The curve $g$ constructed above will in fact be
an approximation of two logarithmic spirals 
(such as those used in the proof of Lemma~\ref{th-logsp}). Then we prove
Theorem~\ref{difexample} in the same way, each time putting 
two rescaled ``hats'' on top of the previous ``hat''. The curve
obtained in this way will not be metrically differentiable at the points
of the Cantor set, since if we consider a sequence of
isosceles triangles $A_nB_nC_n$ with vertex angle
$\angle B_n$ tending to $0$, the ratio between $\|A_nC_n\|$ and 
$\|A_nB_n\|+\|B_nC_n\|$ will
tend to zero as $n\to\infty$, for any norm $\|\cdot\|$. 
\end{remark}

\begin{remark}\label{explanation-next}
If we work with the $\ell_1$-norm, then for a fixed $\alpha\in(0,\pi/2)$ let
$h=\frac{1}{2}\tan\alpha$ and
$$
g(t)=
\begin{cases}
(t+h)-hi, &  \text{ if }t<-h,\\
ti, & \text{ if }t\in[-h,0],\\
t, & \text{ if }t\in[0,1],\\
1-(t-1)i, & \text{ if }t\in[1,1+h],\\
(t-h)-hi, & \text{ if }t>1+h.
\end{cases}
$$

The curve $g$ satisfies conditions of Lemma~\ref{th-logsp}
with any $q<1$ (for the $\ell_1$ norm), 
and although it cannot be made into a graph of a
function in the usual sense,  one can easily see that putting
together such ``boxes'' (rescaling as necessary and taking 
$\alpha_n\to\pi/2$), we obtain the example of a planar curve
with metric derivative $1$ at every point, but with uncountable
set of points where it is not metrically differentiable.
\end{remark}

\section{Metric regularity and metric differentiability}\label{auxsec}

This section contains mainly auxiliary results.
Let $f:[a,b]\to X$, $I=[a,b]$.
We say that $x\in I$ is {\em bilaterally metrically 
regular point of the function~$f$}, provided
\[ \lim_{\substack{(y,z)\to(x,x)\\a\leq y\leq x\leq z\leq b}} 
\frac{\|f(y)-f(z)\|}{\bigvee^{z}_{y} f}=1.\]  
See the beginning of section~\ref{ptsmetrnondifsec} for
the definition of a metrically regular point. Note that
every bilaterally metrically regular point of a function is also its metrically
regular point.

\begin{lemma}\label{mdifandctstomreglem} 
Let $X$ be a Banach space, $g:[a,b]\to X$, $x\in[a,b]$, 
$g$ is metrically differentiable at $x$ with $md(g,x)>0$, 
and $md(g,\cdot)$ is continuous at $x$. Then $x$ is bilaterally metrically 
regular point of the function~$g$.
\end{lemma}

\begin{proof}
Lemma~\ref{ctstomreglem} implies that $x$ is a metrically
regular point of $g$.
Let $\e>0$. By metric differentiability of $g$ at $x$, by Lemma~\ref{kirfedlem}, 
and by continuity of $md(g,\cdot)$ at~$x$ find 
$\delta>0$ such that 
$(1-\e)\,md(g,x)|z-y|\leq\|g(z)-g(y)\|$,
$\bigvee^z_y g=\int^z_y md(g,s)\,ds$,
for $x-\delta<y<x<z<x+\delta$,
and 
$md(g,x+t)<(1+\e)\cdot md(g,x)$
for $|t|<\delta$ with $x+t\in[a,b]$.
Thus, for $y,z$ with $x-\delta<y\leq x\leq z<x+\delta$
we have
\begin{equation*}
\begin{split}
\bigg(\frac{1-\e}{1+\e}\bigg) \bigvee^{z}_{y}g &=	
\bigg(\frac{1-\e}{1+\e}\bigg) \int^{z}_y md(g,s)\,ds
\leq(1-\e)\cdot md(g,x)\,|z-y|\\&\leq\|g(z)-g(y)\|
\leq\bigvee^{z}_{y}g.
\end{split}
\end{equation*}
If $y\neq z$, then by dividing by $\bigvee^{z}_{y}g$, we obtain
$ \frac{1-\e}{1+\e}\leq\frac{\|g(y)-g(z)\|}{\bigvee^{z}_{y}g}\leq 1$,
and thus $x$ is bilaterally metrically regular point of $g$.
\end{proof}

\begin{lemma}\label{ctsplusbiltometrlem}
Let $X$ be a Banach space, and let $f:[a,b]\to X$. If $md(f,\cdot)$
is continuous at $x\in[a,b]$, and $x$ is a~bilaterally metrically
regular point of $f$, then $f$ is metrically differentiable at~$x$.
\end{lemma}

\begin{proof}
If $md(f,x)=0$, then the conclusion follows from Lemma~\ref{mdflem}(i),
and thus we can assume that $md(f,x)>0$.
Lemma~\ref{ctstometric} implies that the condition~\eqref{md} 
holds provided $\sign(z-x)=\sign(y-x)$.
Thus, we only need to treat the case $\sign(z-x)=-\sign(y-x)$
since the cases when either $y=x$ or $z=x$ follow easily from
the existence of $md(f,x)$.
Let $\e>0$. Find $\delta>0$ such that 
for all $x-\delta<y\leq x\leq z<x+\delta$ with $(y,z)\neq(x,x)$ we have that
$\|f(y)-f(z)\| \geq (1-\e)\bigvee^z_y f$,
$\bigvee^z_y f=\int^z_y md(f,t)\,dt$, and 
$(1-\e)\,md(f,t)\leq md(f,x)\leq (1+\e)\,md(f,t)$
for $|x-t|<\delta$ with $t\in[a,b]$.
Let $x-\delta<y\leq x\leq z<x+\delta$.
Then
\begin{equation*}
\begin{split}
\|f(y)-f(z)\| & \geq (1-\e)\bigvee^z_y f =(1-\e)\int^z_y md(f,t)\,dt\\
& \geq (1-\e)^2\,md(f,x)(z-y).
\end{split}
\end{equation*}
The other inequality follows from the same reasoning as in~\eqref{eqreus}.
\end{proof}

\begin{lemma}\label{todiflem}
Let $X$ be a Banach space, let $f:[a,b]\to X$ be continuous, BV, and
such that it is not constant on any subinterval of $[a,b]$. Let $x\in(a,b)$, $y=v_f(x)$, and $g=f\circ v^{-1}_f$.
Then 
\begin{enumerate}
	\item if $x$ is a metrically regular point of $f$,
	then $md(g,y)=1$,
	\item if $x$ is a bilaterally metrically regular point of $f$,
	and there exists a neighbourhood $U$ of $x$ such that all $z\in U$
	are metrically regular points of $f$,
	then $g$ is metrically differentiable at $y$.
\end{enumerate}
\end{lemma}

\begin{proof}
To prove~(i), note that
\begin{equation*}
\label{regudif}
\begin{split}
	1&=\lim_{z\to x} \frac{\|f(z)-f(x)\|}{\big|\bigvee^z_x f\big|}
	 =\lim_{z\to x} \frac{\|f(z)-f(x)\|}{|v_f(z)-v_f(x)|}\\
	 &=\lim_{w\to y} \frac{\|f\circ v_f^{-1}(w)-f\circ v_f^{-1}(y)\|}{|w-y|}
	 =md(g,y).
\end{split}
\end{equation*}
\par
For~(ii), first note that $md(g,y)=1$ by part~(i). 
Let $U$ be the neighbourhood of $x$ such that all $z\in U$
are metrically regular points of $f$. Then part~(i) implies that
$md(g,w)=1$ for all $w=v_f(z)$, where $z\in U$. 
To apply Lemma~\ref{ctsplusbiltometrlem},
it is enough to show that $y$ is a bilaterally metrically
regular point of $g$, but
\begin{equation*}
\begin{split}
\lim_{\substack{(s,t)\to(y,y)\\0\leq s\leq y\leq t\leq v_f(b)}} \frac{\|g(t)-g(s)\|}{\bigvee^t_s g}
&=\lim_{\substack{(s,t)\to(y,y)\\0\leq s\leq y\leq t\leq v_f(b)}} \frac{\|f\circ v^{-1}_f(t)-f\circ v^{-1}_f(s)\|}{t-s}\\
&=\lim_{\substack{(u,v)\to(x,x)\\a\leq u\leq x\leq v\leq b}} \frac{\|f(v)-f(u)\|}{v_f(v)-v_f(u)}=1,
\end{split}
\end{equation*}
where the last equality follows from the fact that $x$ is a bilaterally 
metrically regular point of $f$, and $v_f(v)-v_f(u)=\bigvee^v_u f$ for
any $u,v\in U$, $u<v$ by Lemma~\ref{kirfedlem}.
Now, application of Lemma~\ref{ctsplusbiltometrlem} yields the conclusion.
\end{proof}

We will also need the following simple lemma.

\begin{lemma}\label{forces0lem} 
Let $f:[a,b]\to X$, $x\in[a,b]$, be such that $md(f,x)$ exists,
but $f$ is not metrically differentiable at $x$. Then
if $h$ is a homeomorphism of $[a,b]$ onto itself such
that $f\circ h$ is metrically differentiable at $h^{-1}(x)$,
then $md(f\circ h,h^{-1}(x))=0$.
\end{lemma}

\begin{proof} Lemma~\ref{mdflem} shows that $md(f,x)>0$
(otherwise we have a contradiction with the fact that $f$
is not metrically differentiable at $x$).
Suppose that $h$ is an (increasing) homeomorphism such
that $f\circ h$ is metrically differentiable at $y=h^{-1}(x)$.
For a contradiction, suppose that $md(f\circ h,y)>0$.
Note that
\[ \frac{|h(y+t)-h(y)|}{|t|}=\frac{|h(y+t)-h(y)|}{\|f(h(y+t))-f(h(y))\|}\cdot
\frac{\|f(h(y+t))-f(h(y))\|}{|t|},\]
and it follows that
$ h'(y)=\frac{md(f\circ h,y)}{md(f,x)}>0$.
Thus $h'(y)$ exists and is non-zero. This implies that $(h^{-1})'(x)$
exists. Because $f=(f\circ h)\circ h^{-1}$, Lemma~\ref{mdflem}
implies that $f$ is metrically differentiable at $x$, a contradiction.
We conclude that $md(f\circ h,y)=0$.
\end{proof}

\section{Continuous metric differentiability via homeomorphisms}\label{contsec}

Let $f:[a,b]\to X$.
Let  $M_f$ be the set of all points $x\in[a,b]$ with the following
property:
there is no neighbourhood $U=(x-\delta,x+\delta)$
of $x$  such that either $f|_U$
is constant or all points of $U$ are metrically regular
points of the function $f$.
Obviously, $M_f$ is closed, and $a,b\in M_f$.

\begin{theorem}\label{ctsmdthm}
Let $X$ be Banach space, and let $f:[a,b]\to X$. 
Then the following are equivalent.
\begin{enumerate}
\item There exists a homeomorphism $k$ of $[a,b]$
onto itself such that
$md(f\circ k,\cdot)$ is continuous on $[a,b]$.
\item $f$ is continuous, BV, and $\Hau^1(f(M_f))=0$.
\end{enumerate}
\end{theorem}

\begin{proof} 
To prove that (i)$\implies$(ii), note that 
the existence of continuous metric derivative implies
continuity and boundedness of variation of the function, 
and these properties are preserved when the function
is composed with a homeomorphism.
Thus, it is enough to prove that $\Hau^1(f(M_f))=0$. Note that
$M_f={k}(M_{f\circ {k}})$, and thus it is enough to prove that
$\Hau^1((f\circ {k})(M_{f\circ {k}}))=0$. Let $g=f\circ {k}$.
We claim that 
\begin{equation}\label{Mgeq}
M_g\subset\{x\in[a,b]:md(g,x)=0\}.
\end{equation}
Indeed, Lemma~\ref{ctstomreglem} implies
that every point $x\in(a,b)$, such that $md(g,x)>0$, is metrically regular
point of $g$.
By continuity of $md(g,\cdot)$,
there exists a neighbourhood $U$ of $x$ such that $md(g,y)>0$ 
at all $y\in U$, and thus 
all points of $U$ are metrically regular points of $g$.
So we get~\eqref{Mgeq}, and then by
Lemma~\ref{mdsardlem}, we see that $\Hau^1(g(M_g))=0$.
\par
To prove that (ii)$\implies$(i), let $(U_i)_i$ be the collection
of all maximal open intervals inside
$[a,b]$ such that $f|_{U_i}$
is constant, and put $U=\bigcup_i U_i$.
Define $\vp(t)=v_f(t)+\lambda(U\cap[a,t])$ for $t\in[a,b]$.
Let $(a_j,b_j)$ be the maximal open components of $[a,b]$ 
such that all points of $(a_j,b_j)$ are metrically regular points of $f$.
Let $\alpha_j=\vp(a_j)$, $\beta_j=\vp(b_j)$.
Then $\vp(b_j)-\vp(a_j)=\bigvee^{b_j}_{a_j}f$.
Note that
\begin{equation}\label{ctsmdeq}
\vp(b)=\lambda(U)+\bigvee^b_a f=\lambda(U)+\sum_j\bigvee^{b_j}_{a_j}f=
\lambda(U)+\sum_j(\beta_j-\alpha_j)
=\lambda(\vp([a,b]\setminus M_f)),
\end{equation}
by~\cite[Lemma~2.7]{DZ2}, and thus 
$\lambda(\vp(M_f))=\lambda(M_{f\circ\vp^{-1}})=0$
(the left-hand side of \eqref{ctsmdeq}, 
$\vp(b)$, is equal to $\lambda(\vp[a,b])$, and $\vp$ is increasing).
Let $g=f\circ\vp^{-1}$.
It is easy to see that $g$ is Lipschitz
(because $\vp$ is a homeomorphism).
By Zahorski's lemma (see e.g.\ \cite[p.~27]{GNW}) there exists a continuously differentiable
homeomorphism $h$ of 
$[0,\vp(b)]$ onto itself such that $h'(x)=0$ if and only if $x\in h^{-1}(M_g)$.
Now, by the equality
\begin{equation}\label{splitdif} 
\frac{ g(h(x+t))-g(h(x))}{t}=\frac{g(h(x+t))-g(h(x))}{h(x+t)-h(x)}\cdot
\frac{h(x+t	)-h(x)}{t},
\end{equation}
and by Lemma~\ref{todiflem}, we obtain that $md(g\circ h,x)$ exists
and is continuous at all $x\in\vp(U)\cup\bigcup_j(\alpha_j,\beta_j)$.
By~\eqref{splitdif}, by the choice of $h$ and the fact that $g$ is
Lipschitz,
we easily obtain that $md(g\circ h,x)=0$ for all $x\in h^{-1}(M_g)$,
and that $md(g\circ h,\cdot)=md(f\circ k,\cdot)$ 
is continuous at all such points (where $k=\vp^{-1}\circ h$).
\end{proof}



Let  $M^\mathrm{b}_f$ be the set of all points $x\in[a,b]$ with the following
property:
there is no neighbourhood {$U=(x-\delta,x+\delta)$ }
of $x$  such that either $f|_U$
is constant or all points of $U$ are bilaterally metrically regular
points of the function $f$.
Obviously, $M^\mathrm{b}_f$ is closed and $a,b\in M^\mathrm{b}_f$.

\begin{theorem}\label{mdifctsthm}
Let $X$ be Banach space, and let $f:[a,b]\to X$. 
Then the following are equivalent.
\begin{enumerate}
\item There exists a homeomorphism $h$ of $[a,b]$ onto itself such that
$f\circ h$ is metrically differentiable
at every point of $[a,b]$, and
$md(f\circ h,\cdot)$ is continuous.
\item $f$ is continuous, BV, and $\Hau^1(f(M^\mathrm{b}_f))=0$.
\end{enumerate}
\end{theorem}

\begin{proof} The proof is similar to the proof of Theorem~\ref{ctsmdthm},
and thus we omit it. It uses Lemmas~\ref{mdifandctstomreglem}
and~\ref{todiflem}(ii).
\end{proof}

The following example shows that the scopes of Theorems~\ref{ctsmdthm} and~\ref{mdifctsthm}
are different (see also Remark~\ref{explrem}).

\begin{example} There exists $1$-Lipschitz mapping
$f:[0,1]\to\ell_2$
such that $md(f,x)=1$ for all $x\in[0,1]$, but $f$ is not metrically
differentiable at a dense subset $S$ of $[0,1]$.
\end{example}

\begin{proof} Choose $t_n>0$ with $\sum_n t_n^2=1$,
and $q_n\in(0,1)$ such that $S=\{q_n:n\in\N\}$ is 
dense in $[0,1]$.
Let $f_n:[0,1]\to\R^2$ be defined as
\begin{equation*}
f_n(t)=
\begin{cases}
(t,0)&\text{ for }0\leq t\leq q_n,\\
\frac{(t-q_n)}{\sqrt{2}}\cdot(1,1)+(q_n,0)&\text{ for }q_n<t\leq 1.
\end{cases}
\end{equation*}
It is easy to see that $f_n(0)=0$ and $f_n$ is $1$-Lipschitz for each $n\in\N$.
Define $f:[0,1]\to\ell_2=\sum\oplus_{\ell_2}\ell^2_2$ as $f(t)=(t_n\cdot f_n(t))_n$.
It is easy to see that $f$ is well defined, and $1$-Lipschitz.
First, we will show that $md(f,x)=1$ for all $x\in[0,1]$.
Choose $x\in[0,1]$ and $\e>0$. Find $n_0\in\N$ such that 
$\sum_{n\geq n_0} t^2_n<\e^2$. Find $\delta>0$ such
that $(x-\delta,x+\delta)\cap\{q_j:j\leq n_0\}\subset\{x\}$.
Let $y\in(x-\delta,x+\delta)$ 
and notice that
\begin{equation*} 
\begin{split}
|y-x|&\geq \|f(y)-f(x)\|\\
&=\bigg( \sum_{n\leq n_0} t^2_n\,\|f_n(y)-f_n(x)\|^2_{\ell^2_2} 
+\sum_{n> n_0}t^2_n\,\|f_n(y)-f_n(x)\|^2_{\ell^2_2}\bigg)^{1/2}\\
& \geq \bigg( \bigg(\sum_{n\leq n_0} t^2_n\bigg)^{1/2}-\e\bigg)|y-x|\geq (1-2\e)|y-x|.\\
\end{split}
\end{equation*}
Conclude by sending $\e$ to $0$.
\par
Now we will show that $f$ is not metrically differentiable at
any $x\in S$. Fix $x=q_m\in S$ for some $m$, 
and let $\delta>0$ be such that $0\leq x-\delta<x+\delta\leq 1$.
Then 
\begin{equation*}
\begin{split} 
\frac{\|f(x-\delta)-f(x+\delta)\|}{2\delta}&
=\frac{1}{2\delta}\bigg(t_m^2\|f_m(x-\delta)-f_m(x+\delta)\|^2_{\ell^2_2}\\
&\qquad +\sum_{n\neq m} t^2_n\, \|f_n(x+\delta)-f_n(x-\delta)\|^2_{\ell^2_2}\bigg)^{1/2}\\
&\leq\frac{1}{2\delta}
\bigg(t_m^2 \delta^2 (2+\sqrt{2})+\sum_{n\neq m} 4\delta^2 t^2_n\bigg)^{1/2}\\
&=\bigg( \frac{2+\sqrt{2}}{4}\,t^2_m+\sum_{n\neq m} t^2_n\bigg)^{1/2}=C_m<1,\\
\end{split}
\end{equation*}
and thus $f$ is not metrically differentiable at $x$, as the condition~\eqref{md}
is violated.
\end{proof}

\begin{remark}\label{explrem} 
Lemma~\ref{forces0lem} implies that
if $h$ is a homeomorphism of $[0,1]$ onto itself
such that $f\circ h$ is metrically differentiable
at all $x\in[0,1]$, then $md(f\circ h,y)=0$ for all $y\in h^{-1}(S)$,
which is a dense subset of $[0,1]$. If $h$ could be chosen to further make
$md(f\circ h,\cdot)$ continuous, then $f$ would have to be constant.
Thus, there exists no homeomorphism $h$ of $[0,1]$ onto itself 
such that $f\circ h$
is metrically differentiable at all points of $[0,1]$ while 
$md(f\circ h,\cdot)$ is continuous.
\end{remark}

\section{Acknowledgment} The authors would like to thank 
David Preiss and Lud\v{e}k Zaj\'\i\v{c}ek for many valuable discussions.


\begin{thebibliography}{WWW}



\bibitem[A]{Am}
L.~Ambrosio, \textit{Metric space valued functions of bounded variation}, 
Ann. Scuola Norm. Sup. Pisa Cl. Sci. (4) \textbf{17} (1990), 439--478.

\bibitem[AKh]{AK}
L.~Ambrosio, B.~Kirchheim,
\textit{Rectifiable sets in metric and Banach spaces}, Math. Ann. \textbf{318} (2000), 527--555. 

\bibitem[BL]{BL}
Y. Benyamini, J. Lindenstrauss, {\em Geometric Nonlinear Functional
Analysis, Vol.~1,} Colloquium Publications \textbf{48},
American Mathematical Society, Providence, 2000.

\bibitem[Br]{Br}
A.~Bruckner, {\em Differentiation of real functions}, 
Second edition, CRM Monograph Series, 5, American Mathematical Society, Providence, RI, 1994. 

\bibitem[DP]{DP}
G.~De Cecco, G.~Palmieri, \textit{LIP manifolds: from metric to Finslerian structure},
Math. Z. \textbf{218} (1995), 223--237. 

\bibitem[D1]{Dac}
J.~Duda, {\em Absolutely continuous functions with values in metric spaces},
submitted.

\bibitem[D2]{Dkir}
J.~Duda, {\em Metric and $w^*$-differentiability of pointwise Lipschitz
mappings}, to appear in the Journal of Analysis and its Applications.

\bibitem[DZ]{DZ2}
J.~Duda, L.~Zaj\'\i\v{c}ek, {\em Curves in Banach spaces -- differentiability via
homeomorphisms}, to appear in the Rocky Mountain J. of Math.

\bibitem[EH]{EH}
M.~J.~Evans, P.~D.~Humke, {\em The equality of unilateral derivatives},
Proc. Amer. Math. Soc. \textbf{79} (1980), no.~4, 609--613.

\bibitem[F]{F}
H.~Federer, {\em Geometric Measure Theory,} Grundlehren der math. Wiss., vol.~153, Springer, New York, 1969.

\bibitem[GNW]{GNW}
C.~Goffman, T.~Nishiura, D.~Waterman, {\em
Homeomorphisms in analysis,} Mathematical Surveys and Monographs, vol.~54, 
AMS, Providence, RI, 1997.

\bibitem[J]{Jef}
R.~Jeffery, {\em The theory of functions of a real variable}, Mathematical
Expositions No.~6, University of Toronto Press, 1953.

\bibitem[Kh]{K}
B.~Kirchheim, {\em Rectifiable metric spaces: local structure and
regularity of the Hausdorff measure}, Proc. Amer. Math. Soc.~\textbf{121} 
(1994), 113--123.

\bibitem[Kb]{Ki}
M.~D.~Kirszbraun, {\em Uber die zusammenziehenden und Lipschitzchen
Transformationen}, Fund. Math. \textbf{22} (1934), 77--108.

\bibitem[KS]{KS}
N.~J.~Korevaar, R.~M.~Schoen, 
\textit{Sobolev spaces and harmonic maps for metric space targets},
Comm. Anal. Geom. \textbf{1} (1993), 561--659. 


\bibitem[Z]{Zsur}
L.~Zaj\'\i\v{c}ek, {\em On $\sigma$-porous sets in abstract spaces},
Abstract and Applied Analysis \textbf{2005:5} (2005), 509--534.

\end{thebibliography}
\end{document}